\documentclass[12pt]{article}
\usepackage{amsfonts,amsthm,amsmath,amssymb,amscd,graphicx}
\usepackage{epsfig}
\usepackage{color}

\theoremstyle{plain}
\newtheorem{lemma}{Lemma}

\newtheorem{theorem}{Theorem}

\newtheorem{lm}{Lemma}

\theoremstyle{definition}
\newtheorem{df}{Definition}

\newtheorem{remark}[lemma]{Remark}

\begin{document}


\title{{\Large On bifurcations of area-preserving and non-orientable maps with quadratic homoclinic tangencies }
~\footnote{ This work has been partially supported by the
Russian Scientific Foundation Grant 14-41-00044. Section~3 is carried out by the RSF-grant (project No.14-12-00811). The author~SG has been partially supported  by the by grants of RFBR No.13-01-00589,
13-01-97028-povolzhye and 14-01-00344.
The authors~AD and~MG
has been have been partially supported by the Spanish MINECO-FEDER
Grant MTM2012-31714 and the Catalan Grant 2014SGR504. The author~MG has  been supported by the DFG~Collaborative Research Center TRR~109 ``Discretization in
Geometry and Dynamics''.}}

\author{A. Delshams$^{\dag}$, M. Gonchenko$^\ddag$ and S.V. Gonchenko$^\S$  \\
{\small $^\dag$ Departament de Matem\`atica Aplicada I, Universitat
Polit\`ecnica de Catalunya, Barcelona,
 Spain}\\
{\small $^\ddag$ Institut f\"ur Mathematik, Technische Universit\"at Berlin, Germany}\\
{\small $^\S$  Institute of Applied Mathematics and Cybernetics, Nizhny Novgorod
 University, Russia}\\ {\small {\tt amadeu.delshams@upc.edu, \tt  gonchenk@math.tu-berlin.de, \tt gonchenko@pochta.ru}}
 }

\date{}
\maketitle


\begin{abstract}
We study bifurcations of non-orientable area-preserving maps with quadratic homoclinic tangencies. We study the case when the maps
are given on non-orientable two-dimensional
{surfaces}. We consider one
and two parameter general unfoldings and establish results related to the {emergence} of elliptic
periodic orbits.
\end{abstract}


\section{Introduction}

The main goal of the paper is the {study of} bifurcations of
area-preserving maps (APMs) {defined on} non-orientable {surface}s and {possessing} homoclinic
tangencies.
{For} dissipative  systems,
{the study of bifurcations of homoclinic tangencies is quite traditional}
and many results obtained here have
a fundamental value for the theory of dynamical chaos. One of such
results, known as {\em theorem on cascade of periodic sinks
(sources)}, goes back {already}
to the paper \cite{GaS73}
of Gavrilov and Shilnikov, {see also \cite{N74,G83}}. Note that in
this paper the general case was considered, {when
the initial two-dimensional diffeomorphism has a saddle fixed
point $O$  with multipliers $\lambda$ and $\gamma$, where $0<|\lambda|<1<|\gamma|$ and the \emph{saddle value}
$\sigma\equiv|\lambda||\gamma|$ is not equal to 1, and the invariant  manifolds of $O$ are quadratically tangent
at the points of some homoclinic orbit.}
In this case,
bifurcations of the homoclinic tangency
lead to the appearance of
asymptotically stable (if $\sigma<1$) or completely unstable (if
$\sigma>1$) periodic orbits. Moreover, in any one parameter
general unfolding such orbits are observed for values of the
parameter belonging to an infinite sequence (cascade) of intervals
that do not mutually intersect and accumulate to the value of the
parameter corresponding to the initial homoclinic tangency.

Concerning related results in the conservative
case, we mention, above all, the well known result of S.~Newhouse,~\cite{N77}, on the {emergence of the} so-called 1-elliptic periodic
orbits (there is only one pair multipliers, $e^{i\varphi}$ and
$e^{-i\varphi}$ with $\varphi\neq 0,\pi$, on the unit circle)
under bifurcations of homoclinic tangencies of multidimensional
symplectic maps. However, the Newhouse theorem from~\cite{N77}
does not give answer whether these 1-elliptic periodic
orbits are
generic.\footnote{ The birth of 2-elliptic generic periodic orbits
was proved in \cite{GST98,GST04} for the case of
four-dimensional symplectic maps with homoclinic tangencies to
saddle-focus fixed points.} This fact is {very}
important for the two-dimensional case where an 1-elliptic periodic
orbit is
elliptic and the genericity means the KAM-stability. The birth of {such} generic elliptic periodic
orbits under homoclinic bifurcations
in symplectic two-dimensional maps was established { by Mora and Romero in~\cite{MR97}.}\footnote{ Note that {an} analogous problem
was considered in~\cite{B87,BSh89} when studying
bifurcations of three-dimensional conservative flows with a
homoclinic loop of a saddle-focus equilibrium.}  {In \cite{MR97} a one parameter family $f_\mu$ of two-dimensional symplectic maps was considered such that (i) $f_0$ has a saddle fixed point $O$ with multipliers $\lambda$ and $\lambda^{-1}$, where $0<|\lambda|<1$, and the invariant manifolds of $O$ have a quadratic tangency at the points of a homoclinic orbit $\Gamma_0$, and (ii) this tangency splits generally when $\mu$ varies. Then, as is shown in \cite{MR97}, there exists an  infinite cascade of intervals $\delta_k$, $k=\bar k, \bar k +1,...$ (where $\bar k$ is a sufficiently large integer) of values of $\mu$ such that the intervals  $\delta_k$ accumulate at $\mu=0$ as $k\to\infty$ and  $f_\mu$ has at $\mu \in\delta_k$ a single-round
elliptic orbit of period $k$. Note that all points of this orbit belong to a small fixed neighborhood of the contour $O\cup\Gamma_0$ and any such point is a fixed one for the corresponding first return map $T_k$. It was shown in \cite{MR97} that the map $T_k$ can be written in some rescaled coordinates in the form
\begin{equation}
\bar x = y , \; \bar y = M_k -x - y^2 +\nu_k(y),
\label{HMsymp}
\end{equation}
where the new coordinates $(x,y)$ and parameter $M_k \sim \lambda^{2k}(\mu - \alpha_k)$ can take arbitrary finite values as $k$ are large and the the coefficient $\alpha_k$ and function $\nu_k(y)$ tend to zero as $k\to\infty$ (the latter, along with their derivatives up to order $(r-2)$ if $f_\mu\in C^r$). Thus, the map $T_k$ is asymptotically close (as $k\to\infty$) to the conservative H\'enon map $\bar x = y , \; \bar y = M_k -x - y^2 $ that has an elliptic fixed point when $-1 < M_k <3$. This point is generic for all such values of $M_k$ except for $M_k = 0$ and $M_k = 1/2$ when the map has a fixed point with multipliers $e^{\pm i\pi/2}$ and $e^{\pm i2\pi/3}$, respectively.
This result gives immediately the intervals $\delta_k$ with border points corresponding to $M_k=3$ and $M_k=-1$, respectively.}

{However, the question on coexistence of single-round periodic orbits of different periods, i.e., whether the intervals $\delta_k$ with different numbers intersect, was not considered in \cite{MR97}. Although, this problem is very important, since its solution is only necessary for construction of the bifurcation diagram. We note also that the H\'enon map is degenerate with respect to bifurcations of the fixed point with multipliers $e^{\pm i\pi/2}$ (the strong resonance 1:4), and this would be strange if the same is valid for the first return map $T_k$. Both these problems were solved in the work \cite{GG09}. However, it was required a development of new technics, in particular, the construction of finite-smooth analogous of the analytical Birkhoff-Moser normal form for a saddle map \cite{M56}. First, it was shown in \cite{GG09} that the intervals  $\delta_k$ can intersect indeed and, moreover, they can be even nested. In the latter case \emph{the phenomenon of global resonance}, discovered by S.Gonchenko 
and Shilnikov in \cite{GS01,GS03}, can be observed when, in particular, the symplectic map $f_0$ can have simultaneously infinitely many single-round elliptic periodic orbits of \emph{all successive periods} $k=\bar k, \bar k +1,...$   Second, conditions of nondegeneracy of the resonance 1:4 in the first return maps were found.}

{In principle, these results were obtained by the way of more accurate calculations of the small terms $\alpha_k$ and $\nu_k(y)$ in the rescaled form of $T_k$. As result, the map of form (\ref{HMsymp}) was deduced with
$M_k = - d^{-1}\lambda^{-2k}\left(\mu - \alpha\lambda^k(1 + O(k\lambda^k))\right) + s_0 + O(k\lambda^k)$ and $\nu_k(y) = s_1 \lambda^k y^3 + \lambda^{2k}O(y^4)$, where  $d\neq 0,\alpha, s_0$ and $s_1$ are some coefficients (invariants of a homoclinic structure). Then we can immediately see that if $\mu=0$ and $\alpha=0$ (a codimension two bifurcation case), then $M = s_0 + O(k\lambda^k)$ and, hence, all maps $T_k$ are ``the same'' -- all of them are close to the same map $\bar x = y , \; \bar y = s_0 -x - y^2 $. In this case, if $-1 <s_0 < 3$, every $T_k$ has an elliptic fixed point (with $\varphi$ close to $\arccos\left(1-\sqrt{1+s_0}\right)$). Note that the cubic term $s_1$ responses for nondegeneracy of the resonance 1:4 (if $s_1 \lambda^k > 0$ the point is a saddle with 8 separatrices, if $s_1 \lambda^k < 0$ the point is of elliptic type, i.e. KAM-stable).}

{In this paper we consider area preserving 
maps defined on a non-orientable  surface~${\cal M}_2$ and
study their homoclinic bifurcations.\footnote{{Notice that these maps cannot be symplectic due to the lack of orientation on~${\cal M}_2$}.}
As in the paper~\cite{GG09}, we construct bifurcation diagrams for
single-round periodic orbits
and study the phenomenon of
``global resonance''.
First of all, we
establish
the theorem on cascade of elliptic periodic orbits,
Theorem~\ref{thcasc-n}. However, we note that, unlike the symplectic case,
the first return maps $T_k$ can be nonorientable. We note that, in the case under consideration, the maps $T_k$ are rescaled to a map asymptotically close (as $k\to\infty$) to the nonorientable conservative H\'enon map $\bar x = y,\; \bar y = M +x -y^2$ (see Lemma~\ref{henmain-n}). Thus, $T_k$ can not have elliptic fixed points. However, $T_k^2$ can have and, therefore, cascades from Theorem~\ref{thcasc-n} relate to double-round elliptic points.
In Theorem~\ref{th03-n} we generalize results of Theorem~\ref{thcasc-n} for two-parameter families $f_{\mu,\alpha}$ with governing parameters $\mu$ and $\alpha$ (see formula (\ref{tau}) for $\alpha$) and deduce the result, Theorem~\ref{th:inf}, about the existence of infinitely many double-round elliptic periodic orbits of all successive even periods beginning from some number.}

{We note also that the problem under consideration, as in \cite{GS01,GS03,GG09}, is related to the Poincar\'e conjecture~\cite{PoinMem} on the density of stable (elliptic for the two-dimensional case) periodic orbits in the phase space of non-integrable Hamiltonian systems. This Poincar\'e problem is wide open.\footnote{{In this connection, we note that for (multidimensional) $C^1$-smooth symplectic
diffeomorphisms given on compact manifolds,
the following properties are generic:
1) hyperbolic periodic points are dense in the phase space, \cite{PR83}; 2) every hyperbolic periodic orbit has a transverse homoclinic point in any neighborhood of any point of the phase space, \cite{T70,T72};
3) if a symplectic diffeomorphism $f$ is not Anosov,  then the 1-elliptic
periodic points
of $f$ are dense in the phase space, \cite{N77}.
It is worth remarking that the above-mentioned $C^1$ generic
properties
can become nontypical if one requires a greater smoothness.
Thus, for example, according to the KAM-theory, elliptic periodic orbits
of $C^r$-smooth two-dimensional symplectic diffeomorphisms are generically
stable at $r \geq 5$ (\cite {R70}), whereas, by property 2),
all periodic elliptic orbits of  generic $C^1$-diffeomorphisms
are unstable.}} Therefore, {the above mentioned phenomenon of global resonance (the coexistence of elliptic periodic orbits of all periods)
can be considered as quite relevant.} There are few explicit criteria for the existence of infinitely many elliptic periodic orbits, and this gives a rare opportunity for the construction of Hamiltonian systems with a given structure of stable modes. From this point of view, it is important that in the case considered in our paper, the global resonance is organized in a quite different way: we give an explicit criterion for the existence of elliptic periodic orbits of \emph{ all  even periods}. It is worth mentioning that related problems on the coexistence of large number of elliptic periodic orbits in Hamiltonian systems are quite popular, see e.g. the paper \cite{Tr07} and the corresponding references in it.}

{Note also that  non-orientable APMs appear naturally when restricting a multidimensional symplectic map onto a two-dimensional non-orientable surface as well as when factorizing orientable symplectic maps by a discrete symmetry group. For example, when a map (not necessarily symplectic one) admits certain symmetries, it can be expressed as some even power of a simpler non-orientable map and, thus, the study of the latter map becomes very important, see e.g. \cite{Lamb,DKT(BCh)}.}

\section{Statement of the problem and main results.} \label{sec:Ch1-stat}

Consider a $C^r$-smooth ($r\geq 3$) APM~$f_0$ {defined} on a {\em non-orientable} {surface} ${\cal M}_2$ and
satisfying the following conditions.

\begin{itemize}
\item[{\bf \textsf{A}.}]  $f_0$ has a saddle fixed point $O$ with multipliers
$\lambda$ and $\lambda^{-1}$, where $|\lambda|<1$.

\item[{\bf \textsf{B}.}] $f_0$ has a homoclinic orbit $\Gamma_0$ {where} the stable and unstable invariant manifolds of the saddle
$O$ have a quadratic tangency.
\end{itemize}

Condition~\textsf{A} means that {there exists} a neighborhood (disk)~$U_0$ of the point~$O$ in which the map~$T_0 = f_0|_{U_0}$
is symplectic, i.e., area-preserving and orientable. The map~$T_0$ is called {the} {\em local map}, it is a saddle map that has the
point~$O$ as a fixed one.
{By
condition~\textsf{B}, the stable $W^s$ and unstable $W^u$ invariant manifolds of $O$ intersect non-transversally} at the points of $\Gamma_0$.
Infinitely many such homoclinic
points are {inside}~$U_0$. We take a pair of these points: $M^+ \in W^s_{loc}$ and $M^- \in W^u_{loc}$. Then a natural number $n_0$
exists such that $M^+ = f_0^{n_0}(M^-)$. Let~$\Pi^+\subset U_0$ and $\Pi^-\subset U_0$ be sufficiently small neighborhoods of the
points $M^+$ and $M^-$, respectively. The map~$T_1 = f_0^{n_0}|_{\Pi^-}:\;\Pi^- \to \Pi^+$ is called the {\em global map}. Evidently,
the map $T_1$ is area-preserving (in the symplectic coordinates on $U_0$). However, emphasizing the non-orientability
of phase space,  we assume that

\begin{itemize}
\item[{\bf \textsf{C}.}] The map $T_0$ is symplectic in $U_0$, whereas the map $T_1$ is area-preserving and non-orientable.
\end{itemize}

Let ${\cal H}_0$ be a (codimension one) bifurcation {manifold} composed
of area-preserving $C^r$-maps  {on}~${\cal M}_2$ close to~$f_0$ and such that every map of~${\cal H}_0$ has a nontransversal
homoclinic orbit close to~$\Gamma_0$. Let $f_{\varepsilon}$ be a
family of area-preserving
$C^r$-maps that contains the
map $f_0$ at~$\varepsilon =0$. We suppose that the family depends
smoothly on parameters $\varepsilon =
(\varepsilon_1,...,\varepsilon_m)$ and satisfies the following
condition.

\begin{itemize}
\item[{\bf \textsf{D}.}] The family $f_{\varepsilon}$ is transverse to
${\cal H}_0$ at $\varepsilon=0$.
\end{itemize}

\begin{figure}[htb]
\centerline{\epsfig{file=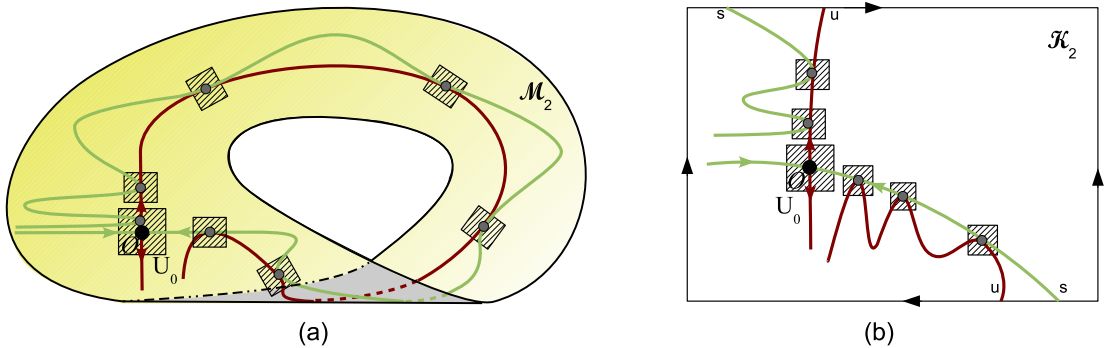,
width=16cm
}} \caption{{\footnotesize An example of non-orientable
APM (a) on the M$\ddot{o}$bius band, (b) on the Klein bottle  having a quadratic homoclinic tangency
at the points of a homoclinic orbit $\Gamma_0$. Some of these
homoclinic points are shown as grey circles. Also a small
neighborhood of the set $O\cup\Gamma_0$ is shown as the union
of a number of ``squares''.}}
\label{fig:fig1}
\end{figure}

Let $U$ be a
small neighborhood of  $O\cup\Gamma_0$ which consists of the small disk $U_0$ containing $O$ and a number of
small disks surrounding those points of $\Gamma_0$ that do not lie
in $U_0$ (see Figure~\ref{fig:fig1}).

\begin{df}
{\em A periodic or homoclinic orbit entirely lying in $U$ is called {\sl p-round}  if it has
exactly $p$ intersection points with any disk of the set $U\backslash U_0$.}
\label{definition:p-round}
\end{df}

In {this paper} we study bifurcations of  {\em
single-round
($p=1$)
periodic orbits} in the families $f_\varepsilon$. Note that every
point of such an orbit can be considered as a fixed point of the
corresponding \emph{first return map}. Such a map is usually
constructed as a superposition $T_k= T_1T_0^k$ of two maps {$T_0\equiv T_0(\varepsilon)= f_\varepsilon\bigl|_{U_0}$}
and
$T_1 \equiv T_1(\varepsilon)= f_\varepsilon^{n_0}: \Pi^- \to U_0 $.
Thus, any fixed point of $T_k$ is a
point of a single-round periodic orbit for $f_\varepsilon$ with
period {$k+n_0$}. We will study maps $T_k$ for all sufficiently
large integer $k$.

Condition~\textsf{D} means that one can introduce in~$U_0$ some {\em canonical coordinates}~$(x,y)$ such that~$J(T_0)
\equiv +1$ and $J(T_1)\equiv -1$, where $J(T)$ is the Jacobian of map~$T$.
{Moreover, we can}
introduce
coordinates~$(x,y)$ in such a way that the
local map $T_0$ can be written in one of the so-called  {\em finitely smooth normal forms}
provided by
the following lemma.

\begin{lemma}\label{lem:nfgen}
{\rm \cite{GG09}}
For any given integer $n$ (such that $n < r/2$ or $n$ is arbitrary for $r=\infty$ or
$r=\omega$--the real analytic case), there is a canonical change of coordinates, of class $C^r$ for
$n=1$ or $C^{r-2n}$ for $n\geq 2$, that brings $T_0$ to the following form
\begin{equation}\label{eq:fng1}
\begin{array}{l}
\bar x = \lambda x\left(1+\beta_1\cdot xy +\dots + \beta_n\cdot
(xy)^n\right) + {x^{n+1}y^n\;O\left(|x|+|y|\right)},\;\;\\
\bar y = \lambda^{-1} y\left(1+\hat\beta_1\cdot xy +\dots +
\hat\beta_n\cdot(xy)^n\right) +
{x^ny^{n+1}\;O\left(|x|+|y|\right)}.\;\;
\end{array}
\end{equation}
The smoothness of these coordinate changes with respect to parameters {can be decreased by 2 unities,} 
i.e., {it is} $C^{r-2}$ for $n=1$ or $C^{r-2n-2}$ for $n\geq 2$, respectively.
\end{lemma}

In these
coordinates, the equations of~$W^s_{loc}\cap U_0$ and~$W^u_{loc}
\cap U_0$ are $y=0$ and $x = 0$, respectively. Moreover,
the normal forms (\ref{eq:fng1}) are very suitable for effective
calculation of maps $T_0^k:(x_0,y_0)\rightarrow (x_k,y_k)$ with
sufficiently large integer $k$. {Indeed}, the following result is valid.

\begin{lemma}\label{lem:Tknfgen}
{\rm \cite{GG09}}
Let $T_0$ be given by (\ref{eq:fng1}), then the map $T_0^k$ can be
written, for any integer $k$, {in the so-called \emph{cross-form}}: 
\begin{equation}
\begin{array}{l}
x_k = \lambda^k x_0\cdot R_n^{(k)}(x_0y_k,\varepsilon)  +
\lambda^{(n+1)k}P_n^{(k)}(x_0,y_k,\varepsilon), \\
y_0 = \lambda^{k} y_k\cdot R_n^{(k)}(x_0y_k,\varepsilon)  +
\lambda^{(n+1)k}Q_n^{(k)}(x_0,y_k,\varepsilon),
\end{array}
\label{eq:Tkgen}
\end{equation}
where
\begin{equation}
\begin{array}{l}
R_n^{(k)}\equiv 1 + \tilde\beta_1(k)\lambda^k x_0y_k +\dots +
\tilde\beta_n(k)\lambda^{nk} (x_0y_k)^n,
\end{array}
\label{eq:Rkn}
\end{equation}
$\tilde\beta_i(k)$, $i=1,\dots,n,$ are some  polynomials (of
degree $i$) with respect to $k$ with coefficients depending on
$\beta_1,\dots,\beta_i,$ {e.g. $\tilde\beta_1 = \beta_1 k$, $\tilde\beta_2  = \beta_2 k - \frac{1}{2}\beta_1^2 k^2$, ...,} and the functions $P_n^{(k)},Q_n^{(k)} =
o\left(x_0^ny_k^n\right)$ are uniformly bounded in $k$ along with
all {their} derivatives {with respect to} coordinates up to order either $(r-2)$ for
$n=1$ or $(r-2n-1)$ for $n\geq 2$.
\end{lemma}

\begin{remark}
{\rm 1)} The normal form of the first order ($n=1$) for $T_0$
\begin{equation}\label{eq:nf1}
\begin{array}{l}
\bar x = \lambda x\left(1+\beta_1\cdot xy\right)  + {x^2y \;O\left(|x|+|y|\right)},\;\;\\
\bar y = \lambda^{-1} y\left(1 - \beta_1\cdot xy\right)  + {xy^2\;O\left(|x|+|y|\right)}\;\;
\end{array}
\end{equation}
is well known from~\cite{MR97,GS90} where it was proved the existence of normalizing
$C^{r-1}$-coordinates. The existence of $C^r$-smooth canonical changes of coordinates (which are
$C^{r-2}$-smooth with respect to parameters) bringing a symplectic saddle map to form~(\ref{eq:nf1}) was proved in~\cite{GST07}.

{\rm 2)} Note that form~(\ref{eq:fng1})
can be considered as a finitely smooth
approximation of the analytical Moser normal form 
\begin{equation}
\begin{array}{l}
\bar x = \lambda(\varepsilon)x\cdot B(xy,\varepsilon),\;\; \bar y
= \lambda^{-1}(\varepsilon) y\cdot B^{-1}(xy,\varepsilon),
\end{array}
\label{eq:BMnf}
\end{equation}
taking place for $\lambda>0$ \cite{M56},  where
$B(xy,\varepsilon)= 1+\beta_1\cdot xy +\dots + \beta_n\cdot (xy)^n
+\cdots$.
Since the form~(\ref{eq:BMnf}) is integrable (e.g. it has integral $xy$), one can easily write the corresponding
formula~(\ref{eq:nf1}) for this case, see~\cite{GS97}.
\label{rem:nf1}
\end{remark}

In {the} coordinates of Lemma~\ref{lem:nfgen}, we can write
$M^+=(x^+,0),M^-=(0,y^-)$.  Without loss of generality, we assume
that $x^{+}>0$ and $y^{-}>0$. Let the neighborhoods $\Pi^{+}$ and
$\Pi^{-}$ of the homoclinic points $M^{+}$ and $M^{-}$,
respectively, be sufficiently small such that
$T_0(\Pi^+)\cap\Pi^+=\emptyset$,
$T_0^{-1}(\Pi^-)\cap\Pi^-=\emptyset$. Then, as usual (see e.g.~\cite{GaS73,book}), the map from
$\Pi^+$ into $\Pi^-$ by orbits of $T_0$ is defined, for all
sufficiently small $\varepsilon$,  on the set consisting of
infinitely many strips $\sigma_k^0\equiv \Pi^+\cap T_0^{-k}\Pi^-$,
$k=\bar k,\bar k+1,\dots$. The image of $\sigma_k^0$ under $T_0^k$
is the strip $\sigma_k^1= T_0^k(\sigma_k^0)\equiv \Pi^-\cap
T_0^{k}\Pi^+$. As $k\to\infty$, the strips $\sigma_k^0$ and
$\sigma_k^1$ accumulate on $W^s_{loc}$ and $W^u_{loc}$,
respectively.

We can write  the global map $T_{1}(\varepsilon):\Pi
^{-}\rightarrow \Pi^{+}$ as follows (in the coordinates of
Lemma~\ref{lem:nfgen})
\begin {equation}
\begin {array}{l}
\overline{x} -x^{+}  =  F(x,y-y^{-},\varepsilon),\;\;
\overline{y} = G(x,y-y^{-},\varepsilon),
\end {array}
\label{eq:t1}
\end{equation}
where $F(0)=0,G(0)=0$. Besides, we have that $G_y(0)=0,
G_{yy}(0)=2d\neq 0$ which follows from the fact (condition~\textsf{B})
that at $\varepsilon=0$
the curve $T_1(W^u_{loc}):\{\overline{x} -x^{+} = F(0,y-y^{-},0),
\overline{y} = G(0,y-y^{-},0)\}$ has a quadratic tangency with
$W_{loc}^s:\{\bar y =0\}$ at $M^+$. When {the} parameters~{$\varepsilon$} vary this
tangency can split and, {moreover, by condition~\textsf{C}, we can introduce the
corresponding splitting parameter as $\mu\equiv
G(0,0,\varepsilon)$.
Accordingly, we can write }
\begin {equation}
\begin {array}{rcl}
F(x,y-y^{-},\varepsilon)  &=&
ax + b(y-y^{-}) + e_{20}x^2 + e_{11}x(y-y^{-}) + e_{02}(y-y^{-})^2 +\; h.o.t., \\
G(x,y-y^{-},\varepsilon)  &=& \mu + cx+d(y-y^{-})^{2} +  f_{20}x^2
+ f_{11}x(y-y^{-}) +  f_{30}x^3 \\
&&+ f_{21}x^2(y-y^{-}) + f_{12}x(y-y^{-})^2 + f_{03}(y-y^{-})^3
+\; h.o.t.,
\end {array}
\label{eq:t1ext}
\end{equation}
where the coefficients $a, b, \ldots, f_{03}$ (as well as $x^+$
and $y^-$) depend smoothly on $\varepsilon$.
Note also that
\begin {equation}
J(T_1) =
\mbox{det}\;\;\left(
\begin {array}{l}
F_x\;\; F_y \\
G_x\;\; G_y
\end{array}
\right) \equiv - 1
\label{eq:det1ext}
\end{equation}
since $T_1$ is non-orientable map by condition~\textsf{C}. In particular, we have
\begin {equation}
\begin {array}{l}
\;\;\;\;\;\;\;\;\;\;\;\;bc \;\equiv \; + 1, \\
\tilde R =  \left(2a + 2e_{02}/bd - bf_{11}/d\right) \;\equiv 0\;
\end {array}
\label{eq:det1ext-JR}
\end{equation}

\begin{figure}[htb]
\centerline{\epsfig{file=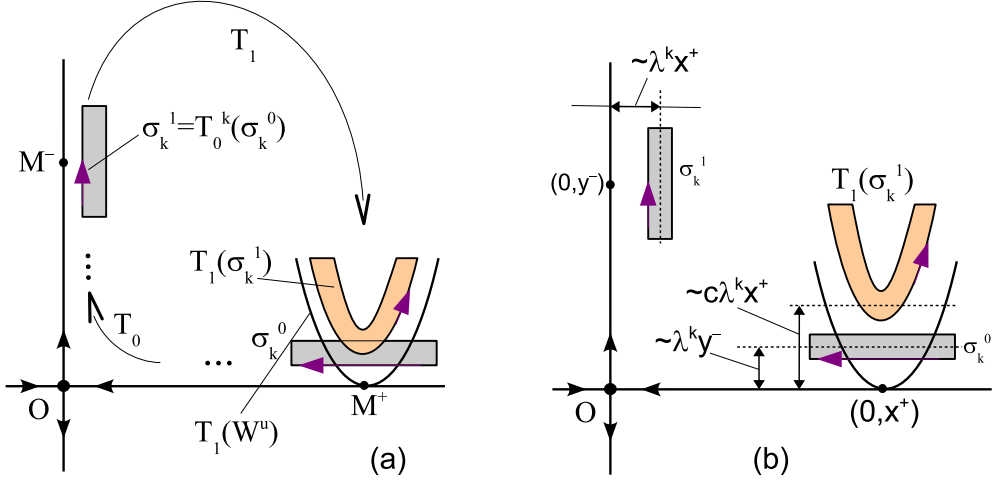, width=15cm }}
\caption{{\footnotesize (a) Actions of the local and
global maps $T_0$ and $T_1$ in $U_0$. Under the map $T_0^k$ some strip $\sigma_k^0$ near $M^+$ is transformed
into a strip $\sigma_k^1$ lying near $M^-$. Under the global map $T_1$ the strip $\sigma_k^1$ is transformed
into a horseshoe $T_1(\sigma_1)$. The latter map is non-orientable, i.e., changes the orientation of the boundary
of the horseshoe with respect to the orientation of the pre-image $\sigma_k^1$. The direction of the orientation
is indicated by the arrows. We assume that $T_0$ is orientable, therefore, the boundary orientation of $\sigma_k^0$
and $\sigma_k^1$ is the same. (b) To a reciprocal position of  the strips
$\sigma_k^0$ and horseshoes $T_1(\sigma_k^1)$. The strips $\sigma_k^0$ and $\sigma_k^1$ are posed on a distance
$\lambda^ky^-(1+...)$ from $W^s_{loc}$ and $\lambda^k x^+(1+...)$ from $W^u_{loc}$, respectively. Hence, the top
of horseshoe $T_1(\sigma_k^1)$ are posed on a distance $c\lambda^k x^+(1+...)$ from $W^s_{loc}$.}}
\label{fig:fretm}
\end{figure}

It is {clear} from~(\ref{eq:t1}) and~(\ref{eq:t1ext}) that
 $\mu$ is
the parameter of splitting of manifolds $W^s(O_\varepsilon)$ and
$W^u(O_\varepsilon)$ with respect to the homoclinic point $M^+$.
Indeed,
the curve $l_u= T_1(W^u_{loc}\cap\Pi^-)$ has the equation
$\displaystyle l_{u}\;:\;\bar y= \mu + \frac{d}{b^{2}}(\bar
x-x^{+})^2(1+O(\bar x-x^{+}))$.
Since the equation of $W^s_{loc}$ is $y=0$ for all (small)
$\varepsilon$, it implies
that the
manifolds $T_1(W^u_{loc})$ and $W^s_{loc}$ do not intersect for $
\mu d>0$, intersect transversally at two points for $\mu d<0$, and
have a quadratic tangency (at $M^+$) for $\mu=0$.
{Besides},
since the strips $\sigma_k^1$ accumulate on the segment
$W^u_{loc}\cap\Pi^-$ as $k\to\infty$, it follows that the images
$T_1(\sigma_k^1)$ of $\sigma_k^1$ under $T_1$ have a horseshoe
form and, moreover, horseshoes $T_1(\sigma_k^1)$ accumulate on
$l_u$ as $k\to\infty$. Therefore, the first return maps $T_k=
T_1T_0^k:\sigma_k^0\to\sigma_k^0$ are, in fact, conservative
horseshoe maps with the Jacobian $-1$. Geometrically, the action of this map looks as in Figure~\ref{fig:fretm}.

When $\mu$ varies near zero infinitely many
bifurcations of horseshoes creation (destruction) occur. In {this} paper
we study these bifurcations and show that they include birth
(disappearance) of {\em elliptic periodic points.}

\begin{figure} [htb]
\centerline{\epsfig{file=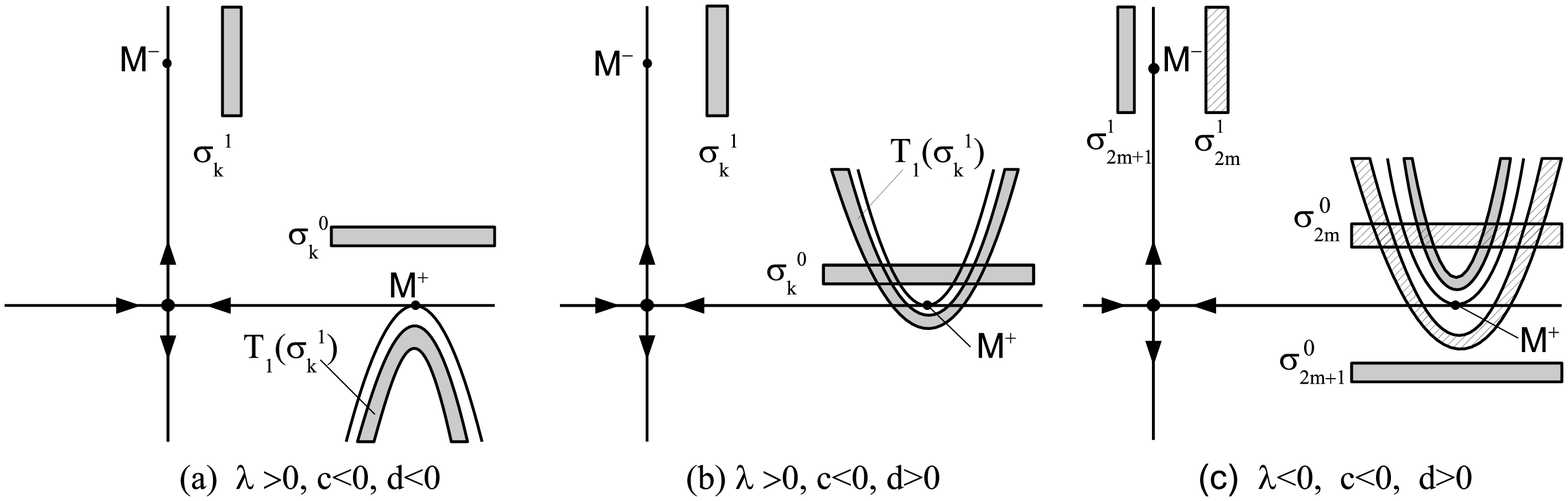, width=16cm
}} \caption{{\footnotesize Types of APMs with a homoclinic
tangency for $c<0$.}} \label{fig:fc<0}
\end{figure}

\begin{figure} [htb]
\centerline{\epsfig{file=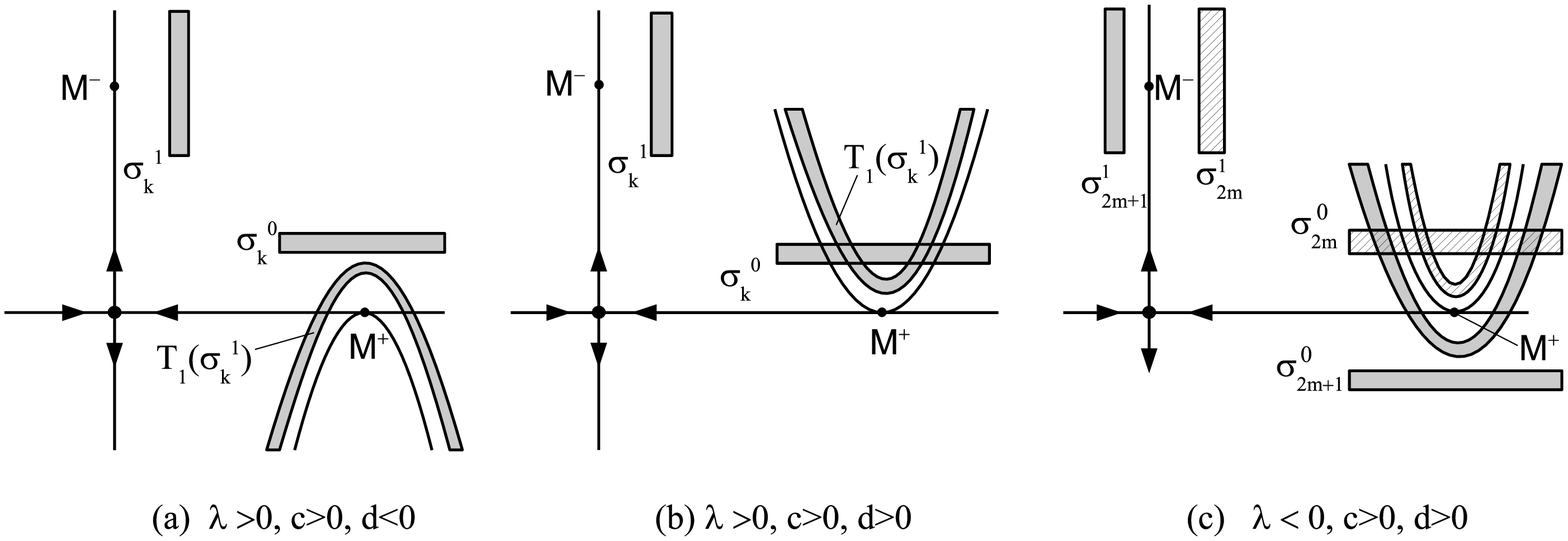, width=16cm
}} \caption{{\footnotesize Types of APMs with a homoclinic
tangency for $c>0$.}} \label{fig:fc>0}
\end{figure}

However, we can also see {that} these horseshoe bifurcations must have
different scenarios depending on {the} type of the initial homoclinic
tangency. Indeed, at $\mu=0$ {the} character of {the} reciprocal position of
the strips $\sigma_k^0$ and their horseshoes $T_1(\sigma_k^1)$ is
essentially defined by the signs of the parameters $\lambda,c$ and
$d$. Moreover, by this feature, we can select 6 different cases of
APMs with quadratic homoclinic tangencies. The
corresponding examples
are shown in
Figures~\ref{fig:fc<0} and \ref{fig:fc>0}. Note that in the cases
with $\lambda<0$ we can always consider $d$ to be positive: if $d$
is negative for the given pair of homoclinic points, $M^+$ and
$M^-$, we can take another pair of points, like $\{T_0(M^+),M^-\}$
or $\{M^+,T_0^{-1}(M^-)\}$, for which the corresponding $d^\prime$
becomes positive.

Note that in the cases with $c<0$, see Figure~\ref{fig:fc<0}, {the}
reciprocal position of all the strips $\sigma_j^0$ and their
horseshoes $T_1(\sigma_j^1)$ at $\mu=0$ is defined quite simply:
$\sigma_j^0\cap T_1(\sigma_j^1)=\emptyset$ if
$\lambda>0,d<0$; the strips $\sigma_j^0$ and  horseshoes
$T_1(\sigma_j^1)$ have regular intersections if
$\lambda>0,d>0$; the corresponding intersections are either
regular for even $j$ or empty for odd $j$ if
$\lambda<0,d>0$. Recall that regular intersection
means here (by \cite{GS87} and \cite{GStT96}) that the set
$\sigma_j^0\cap T_1(\sigma_j^1)$ consists of two connected
components and, moreover, the first return map $T_j\equiv
T_1T_0^j:\sigma_j^0\mapsto\sigma_j^0$ is {\em the Smale horseshoe
map}: its nonwandering set $\Omega_j$ is hyperbolic and
${T_j}\bigl|_{\Omega_j}$ is topologically conjugate to the
Bernoulli shift with two symbols (for more details see
\cite{GStT96,GG09}). Therefore, we can say
that every map $f_0$ in the case $c<0,d>0$ has
infinitely many horseshoes $\Omega_j$, where $j$ runs {for} all
sufficiently large positive integers (respectively, even positive
integers) in the case $\lambda>0$ (respectively, in the case
$\lambda<0$). On the other hand,  every map $f_0$ with $\lambda>0,
c<0, d<0$ has no horseshoes at all (in a small neighborhood $U$).

\begin{figure} [htb]
\centerline{\epsfig{file=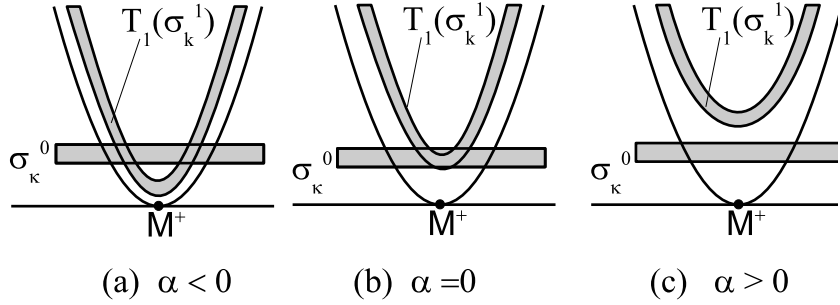, width=12cm
}} \caption{{\footnotesize A horseshoe geometry of symplectic maps
with a homoclinic tangency in the case $\lambda>0,c>0,d>0$ for
various $\alpha$.}} \label{fig:3cltau}
\end{figure}

In the cases of homoclinic tangencies with $c>0$, see
Figure~\ref{fig:fc>0}, {the} reciprocal position of  the strips
$\sigma_j^0$ and horseshoes $T_1(\sigma_j^1)$ depends also on
other invariant {quantities} of the homoclinic structure, {the most important being},
\cite{GS01,GS03,GS87},
\begin{equation}
\label{tau}
\alpha  = \frac{cx^+}{y^-} -1 .
\end{equation}
{First of all, notice that
the sign of $\alpha$ is very important.}
For example, in Figure~\ref{fig:3cltau} it is shown {the}
reciprocal position of the strips $\sigma_j^0$ and horseshoes
$T_1(\sigma_j^1)$ (with sufficiently large $j$) for various values
of $\alpha$ for the case $\lambda>0,c>0,d>0$. Thus, we can see that\footnote{These  results can be easily
explained {by} using the geometry of Figure~\ref{fig:fretm}.
Thus, if $c\lambda^k x^+ > \lambda^k y^-$, i.e., $\alpha>0$, the top of the horseshoe is above the strip
$\sigma_k^0$. {This} means that $f_0$ has no horseshoes and the dynamics is trivial.
However, if $c\lambda^k x^+ < \lambda^k y^-$, i.e., $\alpha<0$, the top of the horseshoe is below the strip
$\sigma_k^0$. Geometrically, it means that $f_0$ has infinitely many horseshoes
(for every sufficiently large $k$). {A} rigorous proof requires quite {elaborate} analytical considerations
which are not presented here, see e.g. \cite{GG09,GS97,GS00,GS03}.}
\begin{itemize}
\item
if $\alpha<0$, then $f_0$
has infinitely many horseshoes $\Omega_j$;
\item
if $\alpha>0$,
then there exists a neighborhood $U(O\cap\Gamma_0)$ in which
{the} dynamics of $f_0$ is trivial: only orbits $O$ and $\Gamma_0$ do
not leave $U$ under iterations of $f_0$.
\end{itemize}
Thus, $\alpha=0$ is {a bifurcation value}, since infinitely many horseshoes appear
(disappear) when varying $\alpha$ near zero (even without splitting
the initial tangency).

Thus, we can draw the following conclusions:
\begin{itemize}
\item[1)]
 the cases of
homoclinic tangencies with $c<0$ or with $\alpha\neq 0$ at $c>0$ are
``ordinary''
and it is sufficient  to study bifurcations of single-round
periodic orbits {only} within {the} framework of one parameter general families
(with the parameter $\mu$);
\item[2)] the cases of homoclinic
tangencies with  $\alpha = 0$ at $c>0$ are ``special''
and it is necessary to
consider at least two parameter general unfoldings (for example,
with parameters $\mu$ and $\alpha$).
\end{itemize}
In {this paper} we adhere to this
approach and present the following {three} theorems as our main results in the case under consideration
(APMs on non-orientable {surface}s).

\begin{figure}[ht]
\psfig{file=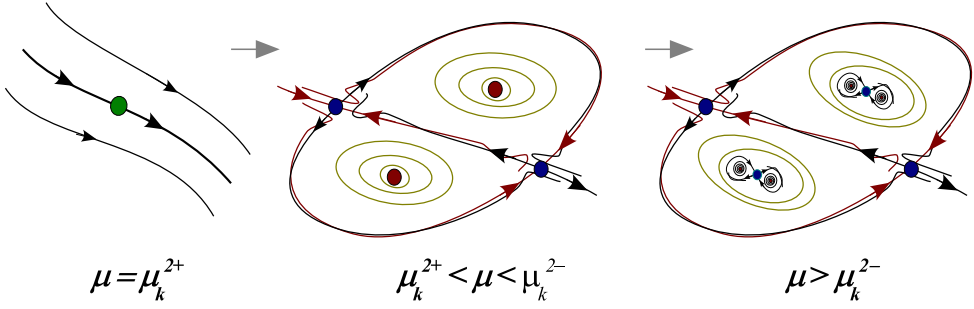, width=16cm}
\caption{{\footnotesize Bifurcation scenarios in the first return maps $T_k$ according to item 3 of
Theorem~\ref{thcasc-n}. We show here that the birth of the elliptic {2-periodic} point {takes place}  when
increasing $\mu$, while for some types of homoclinic tangencies it can occur at decreasing $\mu$.
Since $T_k$ is {a} non-orientable map, then the value $\mu=\mu_k^{2+}$ corresponds to
the appearance of a fixed point with multipliers $+1$ and $-1$. This point {bifurcates} into four points,
two saddle fixed ones and other two points {forming} an elliptic {2-periodic orbit}, when $\mu \in
{\sf e}_k^2$. The {value} $\mu=\mu_k^{2-}$ corresponds to the period doubling bifurcation of this
{elliptic 2-periodic orbit}.
}}
\label{fig:bifscen}
\end{figure}

\begin{theorem} 
Let $f_0$ be an APM satisfying  conditions~\textsf{A}, \textsf{B} and~\textsf{D} and $f_\mu$ be a one parameter family of APMs
that unfolds generally  (under condition~\textsf{C}) at $\mu=0$ the quadratic homoclinic tangency.
{Then for} any interval $(-\mu_0,\mu_0)$ {of} values of $\mu$,
there exists {a positive} integer $\bar k$ {such} that the following holds:

{\rm 1. }
All {APMs} close to $f_0$ have no single-round elliptic periodic orbits, while there exist intervals
${\sf e}_k^2\subset I_\varepsilon$, $k=\bar k,\bar k+1,\dots$,
such that the map $f_\mu$ has at
$\mu\in{\sf e}_k^2$ a
double-round elliptic
{periodic}
orbit, of period $2(k+{n_0})$, which corresponds to a {2-periodic} point of the first return map $T_k$.

{\rm 2.} The intervals ${\sf e}_k^2$ accumulate at $\mu=0$ as $k\to\infty$
and do not intersect for sufficiently large and different integer $k$ if $c<0$, or $\alpha\neq 0$ in the case $c>0$.

{\rm 3.}  Any interval ${\sf e}_k^2$  has border points $\mu=\mu_k^{2+}$
 and $\mu=\mu_k^{2-}$ {where} the map $f_\mu$ has
a single-round periodic  orbit (of period $(k+{n_0})$) with multipliers $+1$ and $-1$ at
$\mu=\mu_k^{2+}$ and
a double-round  periodic  orbit (of period $2(k+{n_0})$) with double multiplier $-1$ at
$\mu=\mu_k^{2-}$. See Figure~\ref{fig:bifscen}.

{\rm 4.} The angular
argument $\varphi$ of the multipliers $e^{\pm i\varphi}$ of the elliptic {periodic orbits for} $\mu\in{\sf e}_k^2$ depends
monotonically on $\mu$ and the elliptic {orbit} is generic (KAM-stable)
for all such $\mu$, except for those {$\varphi=\varphi(\mu)$ such that} $\varphi(\mu)=\frac{\pi}{2},\frac{2\pi}{3},\arccos(-\frac{1}{4})$.
\label{thcasc-n}
\end{theorem}

Note that Theorem~\ref{thcasc-n} does not give answer {to} the question {of the} mutual position of the
intervals ${\sf e}_k^2$ in the critical case $\alpha =0$. But
this {value of $\alpha$} is quite important, since {it is related} to the coexistence of elliptic orbits of different
periods. Therefore, we assume now that $f_0$ is a map satisfying conditions~\textsf{A},~\textsf{B} and~\textsf{D} with $\alpha =0$ and consider a two parameter family
$\{f_{\mu,\alpha}\}$ {which is} a general unfolding for the initial tangency with $\alpha=0$.
Let $D_\epsilon$  be {a} sufficiently small neighborhood (of diameter $\epsilon
>0$) of the origin in the parameter plane $(\mu,\alpha)$. {Then the following result holds.}

\begin{theorem}
{In  $D_\epsilon$, for any}
$\epsilon>0$, {there exist} infinitely many open domains $E_k^2$ (strips)
such that if
%
 $(\mu,\alpha)\in E_k^2$, then  the map $f_{\mu,\alpha}$  has a double-round elliptic orbit of period $2(k+n_0)$ (corresponding
 to a elliptic {2-periodic orbit} of the first return map $T_k$).
The domains $E_{k}^2$ accumulate at the axis $\mu=0$ as $k\to\infty$, all of them are mutually crossed and intersect the axis $\mu=0$.
Every domain $E_k^2$ has two smooth boundaries, {the} bifurcation curves $L_k^{2+}$ and $L_k^{2-}$, which
correspond, respectively, to the existence of a single-round nondegenerate periodic orbit with
double multipliers $+1$  and a double-round nondegenerate periodic orbit with double multipliers
$-1$.
\label{th03-n}
\end{theorem}

\begin{figure}[htb]
\centerline{\epsfig{file=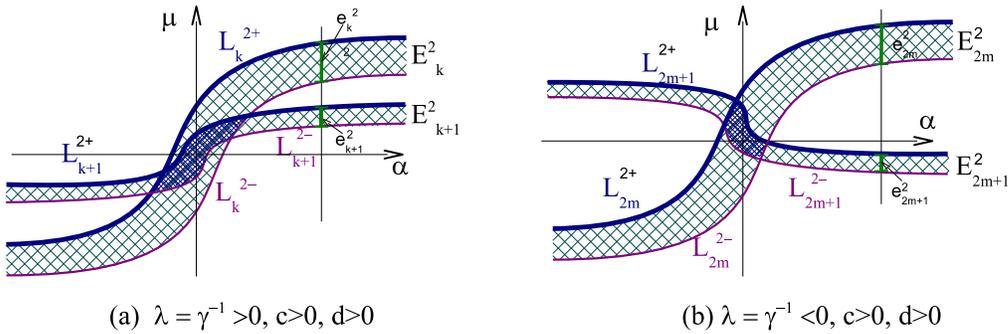, width=16cm
}} \caption{{\footnotesize Elements of the bifurcation diagrams for the families $f_{\mu,\alpha}$.} }
\label{fig2mar-ns-n}
\end{figure}

In Figure~\ref{fig2mar-ns-n} some qualitative illustrations to Theorem~\ref{th03-n} are shown
for the cases where (a) $\lambda>0, c>0, d>0$ and (b) $\lambda<0, c>0, d>0$.

We introduce now the following quantity
\begin{equation}
s_0 = dx^+(ac + f_{20}x^+) - \frac{1}{4}(f_{11}x^+)^2\;,
\label{s0non}
\end{equation}
which is calculated {through the} coefficients of the global map $T_1$, see formula (\ref{eq:t1ext}), and
plays an important role in {the} global dynamics of the map $f_0$ with $\alpha=0$.

Theorem~\ref{th03-n} shows that elliptic (double-round) periodic
orbits of different periods can coexist when values of $\mu$ and $\alpha$ vary near
zero. Moreover, infinitely many such orbits can coexist, in principle, at the global
resonance $\mu=0, \alpha=0$. The following theorem give us sufficient
conditions for this phenomenon.

\begin{theorem} 
Let $f_0$ be {an} APM satisfying  conditions~\textsf{A},~\textsf{B} and~\textsf{D}. We assume also that the resonant condition $\alpha=0$  takes place for
$f_0$. Then, if $s_{0}$ (see formula
(\ref{s0non})) is such that $-1 < s_{0} < 0$, then there exists {a positive} integer $\bar k$ {such} that $f_0$ has
infinitely many double-round
elliptic periodic orbits of all successive even periods $2(k+n_0)$, where $k\geq\bar k$.
Moreover, if  $s_{0}\neq -\frac{1}{2}; -\frac{1}{\sqrt{2}};-\frac{5}{8}$, all these elliptic {periodic} orbits
are generic
(KAM-stable).
\label{th:inf}
\end{theorem}

In the rest part of the paper we prove these and related results.

\section{Rescaling Lemma. }

In principle, for the study of bifurcations in the first return map $T_k= T_1T_0^k$ we could write it in the
initial coordinates using formulas (\ref{eq:t1ext}) and  (\ref{eq:Tkgen})
for the  maps  $T_1$ and $T_0^k$, respectively, and, after, work with the obtained formulas.
However, there is a more effective way for studying bifurcations. Namely, we can bring maps $T_k$
to some unified form for all large $k$ using the so-called rescaling method as it has been done in
many papers.\footnote{see e.g. the papers \cite{BSh89,MR97,GS97,GST98,TR-K98,GS00,GS03} {where}  the
rescaling method was applied for the conservative case.} After this, we can study (one
time) bifurcations in the unified map and ``project'' {the} obtained results onto the first return maps
$T_k$ for various $k$. This ``universal'' map is deduced in the following lemma.

\begin{lm} 
Let $f_\varepsilon$ {be the family} under consideration {satisfying} conditions~\textsf{A--D}. Then, for every sufficiently large $k$,
the first return map $T_k :
\sigma_k^0 \rightarrow \sigma_k^0$ can be brought, by a linear transformation of coordinates and
parameters, to the following form
\begin{equation}
\begin{array}{l}
\bar{X} \; = \; Y, \\
\displaystyle \bar{Y}  \;=\; M +  X - Y^2  + \frac{f_{03}}{d^2}\lambda^{k} Y^3 +
{k\lambda^{2k}\varepsilon_k(Y,M)} \;.
\end{array}
\label{henon0-n}
\end{equation}
{Map (\ref{henon0-n})
is defined} on some asymptotically {large} domain covering in
the limit $k\to+\infty$ all finite values of $X,Y$ and $M$, {the function $\varepsilon_k(Y,M)$} is uniformly
bounded in $k$ {jointly} with all {their} derivatives up to order $(r-4)$), and the following formula takes
place for $M$:
\begin{equation}
\begin{array}{l}
\displaystyle M = -d(1+\nu_k^1)\lambda^{-2k}\left(\mu + \lambda^k(cx^+ - y^-)(1 + k\beta_1
\lambda^k x^+y^-)\right) - s_0 + \nu_k^2\;;
\end{array}
\label{mui-n}
\end{equation}
{with} the coefficient $s_0$ {satisfying} (\ref{s0non}) and {where $\nu_k^{1} = O(\lambda^k), \nu_k^{2} =
O(k\lambda^k)$} are some asymptotically small coefficients.
\label{henmain-n}
\end{lm}

{\em Proof.}
We will use the representation of the symplectic map $T_0$ in the ``second normal form'', i.e., in
form (\ref{eq:fng1}) for $n=2$.\footnote{Clearly, we lose a little in  a smoothness, since
the second order normal form is $C^{r-2}$ only, see Lemma~\ref{lem:nfgen}. However, we get more
important information on form of the first return maps. On the other hand, our
considerations cover also {the} $C^\infty$ and real analytical cases.}  Then the map $T_0^k:\;\sigma_k^0
\rightarrow \sigma_k^1$, for all sufficiently large $k$, can be written in the following form
\begin {equation}
\begin {array}{l}
x_k = \lambda^k x_0(1+\beta_1 k \lambda^k x_0y_k) +
O(k^2\lambda^{3k}),\;\;
y_0 = \lambda^k y_k(1+\beta_1 k \lambda^k x_0y_k) +
O(k^2\lambda^{3k}).
\end {array}
\label{eq:03bnf1}
\end{equation}

Then,
using formulae (\ref{eq:t1ext}) and (\ref{eq:03bnf1}),
we can write the first return  map $T_k : \sigma_k^0 \rightarrow \sigma_k^0$
in the following form
\begin{equation}
\begin{array}{l}
\bar{x}-x^+  =  a\lambda^k x + b(y-y^-) + e_{02}(y-y^-)^2 + \\
\qquad\qquad +O(k|\lambda|^{2k}|x| + |y-y^-|^3 +
|\lambda|^k|x||y-y^-|), \;\; \\ \\
\lambda^{k}\bar{y}\left(1+k\lambda^k\beta_1 \bar{x}\bar{y}\right)
+
k\lambda^{3k} O(|\bar{x}| + |\bar{y}|)\;= \; \\
\qquad\qquad \;=\; \mu + c\lambda^k x\left(1+k\lambda^k\beta_1
xy\right) + d(y-y^-)^2 +
\lambda^{2k}f_{02}\;x^2 + \\
+ \lambda^{k}f_{11}(1+k\lambda^k\beta_1 xy)\;x(y-y^-) +
\lambda^{k}f_{12}\; x (y-y^-)^2 +
f_{03}(y-y^-)^3 + \\
\qquad\qquad + O\left((y-y^-)^4 + \lambda^{2k}|x||y-y^-| +
k|\lambda|^{3k}|x| + k\lambda^{2k}|x||y-y^-|^2\right),
\end{array}
\label{16-n}
\end{equation}
where $x=x_0, y=y_k$.

Below, we will denote by $\alpha_k^i$, $i=1,2,...,$ some {coefficients}
asymptotically small in $k$ such that $\alpha_k^i =
O(k\lambda^k)$. Now we shift the coordinates
$$
\eta = y - y^-,\; \xi = x - x^+ -  \lambda^k x^+(a+\alpha_k^1),
$$
in order to {cancel} the constant term (independent of coordinates) in the first equation of
(\ref{16-n}). Thus, (\ref{16-n}) is recast as follows
\begin{equation}
\begin{array}{l}
\bar{\xi}  = a\lambda^k\xi + b\eta + e_{02}\eta^2 +
O\left(k\lambda^{2k}|\xi| + |\eta|^3 +
|\lambda|^kO(|\xi||\eta|)\right), \\ \\
\lambda^{k}\bar{\eta}(1 + \alpha_k^2)
+ k\lambda^{2k} O(|\bar{\xi}| + \bar{\eta}^2) +
k\lambda^{3k} O(|\bar{\eta}|)
\;=\; \\
\qquad =\;M_1 + c\lambda^k\xi(1  + \alpha_k^3)
+ \eta^2(d + \lambda^k f_{12}x^+) + \lambda^{k}\eta(f_{11}x^+ +
\alpha_k^4) + \lambda^k f_{11}\xi\eta +
%
f_{03}\eta^3  + \\
\qquad\qquad\qquad\qquad + O\left(\eta^4 + k|\lambda|^{3k}|\xi| +
k\lambda^{2k}(\xi^2 + \eta^2) + \lambda^{k}|\xi|\eta^2\right),
\end{array}
\label{sys11-n}
\end{equation}
where
\begin{equation}
M_1 = \mu + \lambda^k(cx^+ - y^-)(1 + k\lambda^k\beta_1 x^+y^-) + \lambda^{2k}x^+(ac + f_{02}x^+) +
O(k\lambda^{3k}).
\label{m1or}
\end{equation}

Now, we rescale the variables:
\begin{equation}
\displaystyle \;\xi = -\frac{b(1+\alpha_k^2)}{d +
\lambda^kf_{12}x^+}\lambda^{k} u \;,\; \eta =
-\frac{1+\alpha_k^2}{d + \lambda^kf_{12}x^+}\lambda^{k} v .
\label{reskk-n}
\end{equation}
System (\ref{sys11-n}) in coordinates $(u,v)$ is rewritten in the
following form
\begin{equation}
\begin{array}{l}
\displaystyle \bar{u}  = v + a\lambda^k u -
\frac{e_{02}}{bd}\lambda^{k} v^2 +
O(k\lambda^{2k}), \\ \\
\bar{v} \; = \;
M_2 +  u(1 + \alpha_k^5) - v^2 + \\
\displaystyle \qquad\qquad + v(f_{11}x^+ + \alpha_k^6)  -
\frac{f_{11}b}{d}\lambda^{k} uv + \frac{f_{03}}{d^2}\lambda^{k}
v^3 + O(k\lambda^{2k}) \;,
\end{array}
\label{0017-n}
\end{equation}
where
$$
\displaystyle M_2 = -\frac{d +
\lambda^kf_{12}x^+}{1+\alpha_k^2}\lambda^{-2k} M_1.
$$
The following shift of coordinates (we remove the {terms} linear in $v$ from the second equation)
$$
u_{new} = u - \frac{1}{2}(f_{11}x^+ + \alpha_{k}^6),\; v_{new} = v
- \frac{1}{2}(f_{11}x^+ + \alpha_{k}^6),
$$
brings map (\ref{0017-n}) to the following form
\begin{equation}
\begin{array}{l}
\displaystyle \bar{u}  = v + a\lambda^k u -
\frac{e_{02}}{bd}\lambda^{k} v^2 +
O(k\lambda^{2k}), \\ \\
\displaystyle \bar{v}  = M_3 + u - v^2  - \frac{f_{11}b}{d}\lambda^{k} uv +
\frac{f_{03}}{d^2}\lambda^{k} v^3 + O(k\lambda^{2k}) \;,
\end{array}
\label{00170-n}
\end{equation}
where
$$
M_3 = M_2  + \frac{(f_{11}x^+)^2}{4}.
$$

Now, we make the following linear change of coordinates
\begin{equation}
x \;=\; u + \tilde\nu_k^1 v \;\;,\;\; y \;=\; v - \tilde\nu_k^2 u
\;, \label{110c}
\end{equation}
where
\begin{equation}
\displaystyle \tilde\nu_k^1 = -\frac{e_{02}}{bd}\lambda^k, \;
\tilde\nu_k^2 = -\frac{e_{02}}{bd}\lambda^k - a\lambda^k.
\label{110d}
\end{equation}
Then, system (\ref{00170-n}) is rewritten as
\begin{equation}
\begin{array}{l}
\bar{x}  =  y + M_3\tilde\nu_k^1 +
O(k\lambda^{2k}), \\ \\
\displaystyle \bar{y}  = M_3 +  x - y^2  + a\lambda^k y - \tilde R\lambda^{k} xy +
\frac{f_{03}}{d^2}\lambda^{k} y^3 + O(k\lambda^{2k}) \;,
\end{array}
\label{0018-n}
\end{equation}
where $ \tilde R =  \left(2a + 2e_{02}/bd - bf_{11}/d\right) \equiv 0 $ by (\ref{eq:det1ext-JR}).
Hence, map (\ref{0018-n}) has the following form
\begin{equation}
\begin{array}{l}
\bar{x}  =  y + M_3 \tilde\nu_k^1 +
O(k\lambda^{2k}), \\ \\
\displaystyle \bar{y}  = M_3 + x - y^2  + a\lambda^k y +
\frac{f_{03}}{d^2}\lambda^{k} y^3 + O(k\lambda^{2k}) \;,
\end{array}
\label{0019}
\end{equation}

Finally, we make one more shift of coordinates
$$
X = x - \frac{1}{2}a\lambda^k - \tilde\nu_k^1 M_3,\;\; Y = y -
\frac{1}{2}a\lambda^k,
$$
in order to {cancel} in (\ref{0019}) the constant term in the first equation and the {term} linear in $y$
 in the second equation. After this, we obtain the final form (\ref{henon0-n}) of map $T_k$ in
the {rescaled} coordinates where formula (\ref{mui-n}) takes place for the parameter $M$. $\;\;\Box$

\section{Proofs of {the} main results.}

{The} bifurcations in {the} first return maps $T_k$ can be studied now {by} using their normal forms deduced
{from} the  rescaling Lemma~\ref{henmain-n}.
Since these normal forms coincide up to asymptotically small {terms} as $k\to\infty$ with the non-orientable
conservative H\'enon map, we recall in the next section some necessary results on bifurcations of fixed points in one parameter
families of conservative H\'enon map in {the} non-orientable  case.

\subsection{On bifurcations of fixed points in the conservative non-orientable H\'enon maps.}

{The} rescaling Lemma~\ref{henmain-n} shows that the unified
limit form for the first return maps $T_k$ is the non-orientable and conservative
H\'enon map
\begin{equation}
\begin{array}{l}
 \bar x = y, \; \bar y = M + x - y^2,
\end{array}
\label{conshen-0}
\end{equation}
with Jacobian $J=-1$.
Bifurcations of fixed points in the conservative H\'enon family
are well known.

Since the H\'enon map (\ref{conshen-0}) is not orientable, it cannot have elliptic fixed
points. However, {there exist} elliptic {2-periodic orbits} for $M\in (0,1)$. The map has no fixed points
for $M<0$, it has one fixed point $\bar O(0,0)$ with multipliers $\nu_1=+1,\nu_2=-1$ at $M=0$ and
two saddle fixed points ($\bar O_1(-\sqrt{M},-\sqrt{M})$ and $\bar O_2(\sqrt{M},\sqrt{M})$) at
$M>0$. Besides, an elliptic {2-periodic} orbit exists for $0<M<1$, {consisting} of the two points $
p_1=(-\sqrt{M},\sqrt{M})$ and $p_2=(\sqrt{M},-\sqrt{M})$; the value $M=1$ corresponds to {a} period
doubling bifurcation of this orbit.   
Note that
the elliptic {2-periodic} orbit is generic for all $M\in (0,1)$ except for $M= \frac{1}{2}$ and {$M=
\frac{3}{4}$} which correspond to the strong resonances $1:4$ and $1:3$, respectively, and $M
= \frac{5}{8}$ which corresponds to {the cancellation of} the first Birkhoff coefficient at the cycle
$\{p_1,p_2\}$, see~\cite{DGGT12}.

It is also known (see, e.g., \cite{DN78,AS82}) that if $M > 5 + 2 \sqrt{5}$ (this is only a
sufficient condition), then the nonwandering set of map (\ref{conshen-0}) is a Smale horseshoe which is
non-orientable for this case.

\subsection{
Proof of Theorem~\ref{thcasc-n}.}\label{sec:fmu}

The proof is deduced from the  rescaling lemma~\ref{henmain-n}. Indeed, since {the} bifurcations of fixed points of
the H\'enon map (\ref{conshen-0})
are known, we can use this information directly {to recover the}
bifurcations of {the} single-round periodic orbits in the family $f_\mu$. We {only} need to know {the} relations
between the parameters of the rescaled map (\ref{henon0-n}) and the initial parameters (i.e., in fact,
between $M$ and $\mu$).

In the case under consideration, the relations between $M$ and $\mu$ are given by formula
(\ref{mui-n}) from which we find $\mu$ as follows
\begin{equation}
\begin{array}{l}
\displaystyle \mu  = - \lambda^k y^- \alpha (1 + k\beta_1 \lambda^k x^+y^-) - \frac{1}{d}(M + s_0 +
\hat\rho_k^1)\lambda^{2k},
\end{array}
\label{mui1-n}
\end{equation}
where $\hat\rho_k^1 = O(k\lambda^k)$ is some small coefficient and
$\displaystyle \alpha = \frac{cx^+}{y^-} - 1$ (see formula (\ref{tau})).

As it follows from Lemma~\ref{henmain-n}, the conservative non-orientable H\'enon map
$\bar{x} \; = \; y, \;\; \bar{y}  \;=\; M + x - y^2, \;$
where $M$ satisfies (\ref{mui-n}), is {the} normal (rescaled) form for the first return maps $T_k$ with all sufficiently large $k$.
This H\'enon map  has no elliptic fixed  points, however, {there exists} an
elliptic {2-periodic orbit} for $0<M<1$. Thus, we obtain, by (\ref{mui-n}), that the first return map
$T_k$ has a fixed point with multipliers $\nu_1=+1,\nu_2 = -1$ (i.e., when $M=0$) if
\begin{equation}
\begin{array}{l}
\displaystyle \mu  = \mu_k^\pm  = - \lambda^k y^- \alpha(1 + k\beta_1 \lambda^k x^+y^-) -
\frac{1}{d}(s_0 + \hat\rho_k)\lambda^{2k},
\end{array}
\label{muk+-1}
\end{equation}
and a {2-periodic orbit} with multipliers $\nu_1=\nu_2 = -1$ (i.e., when $M=1$) if
\begin{equation}
\begin{array}{l}
\displaystyle \mu  = \mu_k^{2-} = - \lambda^k y^- \alpha (1 + k\beta_1 \lambda^k x^+y^-) -
\frac{1}{d}(s_0 + 1 + \hat\rho_k)\lambda^{2k}.
\end{array}
\label{muk--1}
\end{equation}
Thus, the first return map $T_k$ has in this case {an elliptic 2-periodic orbit} when
$\mu\in{\sf e}_k^2$, where ${\sf e}_k^2$ is the interval of values of $\mu$ with border
points $\mu=\mu_k^\pm$ and $\mu=\mu_k^{2-}$. Evidently, if $\alpha\neq 0$,  the intervals ${\sf e}_k^2$ with sufficiently
large and different $k$ do not intersect. $\;\;\Box$

\subsection{
Proofs of Theorems~\ref{th03-n} and \ref{th:inf}.}

{\em Proof of Theorem~\ref{th03-n}.}
By (\ref{muk+-1}) and (\ref{muk--1}), the equations of {the} bifurcation curves $L_k^{2+}$ and $L_k^{2-}$, which are
boundaries of the domain $E^2_k$, can be written as follows
%
\begin{equation}
\begin{array}{l}
\displaystyle L_k^{2+}:\;\;  \mu  = - \lambda^k y^- \left(\frac{cx^+}{y^-} - 1 \right)(1 + k\beta_1 \lambda^k x^+y^-) -
\frac{s_0 + \cdots}{d}\lambda^{2k},
\end{array}
\label{Lmuk+-1}
\end{equation}
%
%
\begin{equation}
\begin{array}{l}
\displaystyle L_k^{2-}:\;\; \mu  - \lambda^k y^- \left(\frac{cx^+}{y^-} - 1 \right) (1 + k\beta_1 \lambda^k x^+y^-) -
\frac{1+ s_0 + \cdots}{d}\lambda^{2k}.
\end{array}
\label{Lmuk--1}
\end{equation}
Since $\lambda^{2k}\ll\lambda^k$, the domains $E^2_k$ with sufficiently large $k$ {do} not mutually
intersect and do not intersect the axis $\mu=0$, if $cx^+\neq y^-$. Thus, the domains do not {always} intersect in the cases with $c<0$.
However, at the global resonance $\alpha = (cx^+/y^- -1) =0$ (which is possible only when $c>0$), as it follows from
(\ref{Lmuk+-1}) and (\ref{Lmuk--1}), all {the} domains $E^2_k$ with sufficiently large $k$ {do}
mutually intersect and all {of} them intersect
the axis $\mu=0$ (as in
Figure~\ref{fig2mar-ns-n}). $\;\;\Box$

{\em Proof of Theorem~\ref{th:inf}.} Assume, for more definiteness, that $d>0$ for all {the} cases under consideration. The case
$d<0$ is treated in the same way.
{Assume that} $f_0$ satisfies conditions~\textsf{A},~\textsf{B} and~\textsf{D} and {that the resonant condition $\alpha =0$ takes place for $f_0$}.
Then, for the one parameter family $f_\mu$ with fixed $\alpha=0$,
 the intervals ${\sf e}_k^2$ have, by (\ref{muk+-1})--(\ref{muk--1}), {the} form
$$
{\sf e}_k^2 = (-1-s_0, - s_0)\frac{\lambda^{2k}}{d}.
$$
Evidently, if $-1<s_0<0$, these intervals will be nested and containing $\mu=0$. {This} implies that the diffeomorphism $f_0$
has infinitely many double-round elliptic periodic orbits.

As  follows from Lemma~\ref{henmain-n}, all the first return maps $T_k$ (with sufficiently large $k$) are reduced to the same
rescaled normal form --- the non-orientable H\'enon map $\bar x = y,\;\bar y = - s_0 + x - y^2$. It is well known that,
for $-1<s_{0}<0$, the {elliptic 2-periodic orbit} of this map is generic if $s_{0}\neq -\frac{1}{2}; -\frac{3}{4};-\frac{5}{8}$.
{These} exceptional cases {are related}, respectively, to {the} resonances $1:4$, $1:3$ and
{to an} elliptic point (at $s_0 = -\frac{5}{8}$) whose
first Birkhoff coefficient is zero.  $\;\;\Box$ 

\subsection*{Acknowledgements.} The authors thank D.Turaev and L.Lerman
for fruitful discussions and remarks.

\end{document}